\documentclass[12pt,twoside]{amsart} \usepackage{amssymb}

\numberwithin{equation}{section}
\newtheorem{theorem}{Theorem}[section] \newtheorem{lemma}[theorem]{Lemma}
\newtheorem{proposition}[theorem]{Proposition}
\newtheorem{corollary}[theorem]{Corollary}

 \theoremstyle{definition}
\newtheorem{definition}[theorem]{Definition}
\newtheorem{remark}[theorem]{Remark}

 \newcommand\codim{\operatorname{codim}}
\newcommand\coker{\operatorname{coker}} \newcommand\Ann{\operatorname{Ann}}
 
\newcommand\Hom{\operatorname{Hom}} \newcommand\Ext{\operatorname{Ext}}

\newcommand{\Proj}{\operatorname{Proj}} 
 \newcommand\rad{\operatorname{rad}}

\newcommand{\cL}{{\mathcal L}}

\newcommand {\PP}{\mathbb{P}}

\begin{document}
\title[Normal form for space curves in a double plane]{Normal form for
space curves
in a double plane}
\author[Nadia Chiarli, Silvio Greco, Uwe Nagel]{Nadia Chiarli$^*$, Silvio
Greco$^*$,
Uwe Nagel$^*$}
\address{Dipartimento di Matematica, Politecnico di Torino,
I-10129 Torino, Italy}
\email{nadia.chiarli@polito.it}
\address{Dipartimento di Matematica, Politecnico di Torino, I-10129 Torino,
Italy}
\email{silvio.greco@polito.it}
\address{Department of Mathematics,
University of Kentucky, 715 Patterson Office Tower,
Lexington, KY 40506-0027
}
\email{uwenagel@ms.uky.edu}

\thanks{$^*$ Partially supported by INDAM-GNSAGA. \\
The authors also gratefully acknowledge partial support by MIUR (N.\ C.,
S.\ G.) and by
a Faculty Summer
  Research Fellowship from the University of Kentucky (U.\ N.), respectively}

\begin{abstract}
This note is an attempt to relate explicitly the geometric and algebraic
properties of
a space curve that is contained in some double plane. We show in particular
that the
minimal generators of the homogeneous ideal of such a curve can be written
in a very
specific form. As applications we characterize the possible Hartshorne-Rao
modules of
curves in a double plane and the minimal curves in their even Liaison classes.
\end{abstract}

\maketitle
\tableofcontents
\section{Introduction} \label{section-intro}
\bigskip
Hartshorne and Schlesinger \cite{HS} have shown that a curve $C
\subset \PP^3_K$ (the projective three-space over an algebraically
closed field $K$) contained in a double plane $2 H$ determines uniquely a
triple of subschemes $Z \subset C' \subset P$, all contained in the
plane $H$, where $C'$ and $P$ are (planar) curves and $Z$ is either empty or
zero-dimensional, and $Z$ is empty if and only if $C$ is arithmetically
Cohen-Macaulay
(aCM).
This correspondence is not bijective.
The goal of this paper is to make the relationship very explicit. In
fact, we show that every non-aCM curve $C$ in a double
plane is related to a certain matrix which allows to describe explicitly the
homogeneous ideal of the curve $C$, as well as to read off the pair
$(Z, C')$ and the Hartshorne-Rao module of
$C$. Moreover, our description of the ideal clarifies the
role of $P$ and allows to compute a graded minimal free resolution
of $C$ (including the maps). Thus, we think of our
description of the ideal as a normal form for curves in a double plane.
In order to explain the relations just mentioned we may assume
that the plane $H$ is defined by the ideal $x R$ where $R = K[x, y,
z, t]$ is the coordinate ring of $\PP^3$. Thus, $Z$ being a 2-codimensional
subscheme
of $H$,  its homogeneous ideal is the ideal of the maximal minors
of its Hilbert-Burch matrix $A \in S^{s, s+1}$ with entries in $S
= K[y, z, t] \cong R/x R$. Denote by $(x, p)$ the homogeneous
ideal of the curve $C'$. Since $Z$ is contained in $C'$ we can write $p$
as the determinant of a homogeneous matrix
\begin{equation*}\label{expression of p}
p = (-1)^s \det \left [\begin{array}{c} A\\ p_1\;\dots \;p_{s+1}
\end{array}\right ]
\end{equation*}
where $p_1,\ldots,p_{s+1} \in S$ are suitable homogeneous polynomials and the
sign is used to simplify the presentation
below. Moreover, the specific form of the ideal of the curve $C$ determines
a vector $^t(f_1,\ldots,f_{s+1})$ such that the matrix
$$
B:= \left[\begin{array}{cc} A&\begin{array}{c}
f_1\\ \vdots
\end{array}\\ p_1\;\dots \;p_{s+1}&f_{s+1}
\end{array}
\right] \in S^{s+1, s+2}
$$
is homogeneous and the last column has in particular maximal degree
among the columns of $B$.
For any $u \times (u+1)$ matrix $N$ we denote by $N_i$ the
determinant of the matrix obtained from $N$ by deleting its $i$-th
column.  This notation allows us to summarize part of our results as
follows:

\begin{theorem} \label{thm-intro} \mbox{}
\begin{itemize}
\item[(a)] Every non-aCM curve $C$ in a double plane determines a
homogeneous matrix
\begin{equation}\label {matrix B}
B:= \left[\begin{array}{cc} A&\begin{array}{c}
f_1\\ \vdots
\end{array}\\ p_1\;\dots \;p_{s+1}&f_{s+1}
\end{array}
\right] \in S^{s+1, s+2}
\end{equation}
where $s \ge 1$, the maximal minors of $A = (a_{i,j})$ generate the ideal
of $Z$ as a
subscheme of
$H$,
$B_{s+2} \neq 0$ defines $C' \subset H$, $\deg f_1 \ge \deg a_{1,1} +
\sum_{j=1}^s
\deg a_{j,j+1}-1$  and
\begin{equation}\label{submatrix M}
M:= \left[\begin{array}{cc} A&\begin{array}{c}
f_1\\ \vdots \\
f_{s}
\end{array}\\
\end{array}
\right] \in S^{s, s+2}
\end{equation}
is a presentation matrix of the Hartshorne-Rao module of $C$ as an $S$-module.
\medskip

\item[(b)] Conversely, let $B \in S^{s+1, s+2}, s \ge 1$ be a
homogeneous matrix as in (\ref {matrix B}) and let $M$ be its
submatrix as in \ref {submatrix M}. Assume:
\begin{itemize}
\item[(i)] the maximal minors of $M$ do not have common zeroes in
$H$;
\item[(ii)] $B_{s+2} \neq 0$;
\item[(iii)] $\deg f_1 \ge \deg a_{1,1} + \sum_{j=1}^s \deg
a_{j,j+1}-1$.
\end{itemize}
Put $p = (-1)^s B_{s+2}$ and let $h \in S$ be a non-zero
homogeneous polynomial with $\deg h = \deg f_1 - (\deg a_{1,1} +
\sum_{j=1}^s \deg
a_{j,j+1}-1)$. Then the ideal
$$
(x^2, x p, p h A_1 + x B_1,\ldots, p h A_{s+1} + x B_{s+1})
$$
defines a curve in the double plane $\{x^2 = 0\}$ associated to
the triple $Z \subset C' \subset P$ where $Z, C'$ and the Hartshorne-Rao module
of $C$ are determined as in part {\rm (a)} and $P$ is defined by $(x, p h)$.
\end{itemize}
\end{theorem}

As already mentioned, writing the generators of the ideal of $C$ as
above allows to compute the minimal free resolution of $C$ (Theorem
\ref{resolution of
J}).  Moreover, it also allows to
describe all possible Hartshorne-Rao modules of curves in a double plane,
as indicated
below.

\begin{corollary} \label{cor-intro}
A graded $R$-module of finite length is the Hartshorne-Rao module of a
curve in the
double plane $\{x^2 = 0\}$ if and only if it can be represented by
a homogeneous matrix
$$
\left [ x E_s \; M \right ] \in R^{s, 2s+2}
$$
where the entries of $M$ are in $(y, z, t) S$, satisfy the degree condition
(iii) of Theorem \ref {thm-intro}~(b) and the maximal minors of $M$
do not have a common zero in $H$ (we denote by $E_s$  the $s\times s$
identity matrix).
\end{corollary}

It is well-known that the even Liaison class of a space curve is determined by
its Hartshorne-Rao module. Since the whole even Liaison class can be
obtained from
a minimal curve by basic double links and possibly a flat deformation, one is
interested in determining the minimal curves. Below, we characterize when a
minimal curve
must lie in some double plane.

\begin{proposition} \label{prop-intro} Let $N \neq 0$ be a graded
$R$-module of
finite length and denote by $\cL_N$ the even liaison class
determined by
$N$.
Then the following conditions are equivalent:
\begin{itemize}
\item[(i)] $\cL_N$ contains a curve lying in some double plane $2 H$;
\item[(ii)] every minimal curve of $\cL_N$ lies in the double
plane $2 H$ or, in case $N$ is cyclic and  annihilated by two
independent linear forms, is an extremal curve;
\item[(iii)] $N$ is one of the modules described in Corollary \ref{cor-intro}.
\end{itemize}
\end{proposition}

The starting point for our results is the observation that every
curve in a double plane sits in the middle of a residual sequence.
We call this sequence the ``expected residual sequence''.
In Section \ref{sec-resid-seq} we study and characterize the
one-dimensional schemes which fit into such an expected residual sequence.
For them
we also determine a graded minimal free resolution. In Section \ref{curves}
we analyze
when such a scheme is indeed a curve, i.e.\ locally Cohen-Macaulay. This
leads to
the normal form for a curve in a double plane. Finally, we discuss the
Liaison class of such a curve.
\bigskip

\section{Standing Notation and Preliminaries} \label{standing}
\bigskip

Throughout the paper we will use the following notation.

\begin{itemize}
\item $K$ is an algebraically closed field.
\item $R: = K[x,y,z,t]$, $S := K[y,z,t] \subseteq R$. We shall identify $S$
with
$R/xR$ when no confusion is possible.
\item $C \subseteq \PP^3 := \Proj(R)$ is a non-degenerate, projective curve of
degree $d$
and arithmetic genus $g$, where by curve we mean a pure $1$-dimensional
locally Cohen-Macaulay projective subscheme (i.e. without
zero-dimensional components).
\item $H$ is the plane defined by the ideal $xR$ and $2H$ denotes the closed
subscheme of
$\PP^3$ whose homogeneous ideal is $(x^2)$. It is called {\it double plane}.
\item If $X \subseteq \PP^3$ is a closed subscheme, then
$I_X \subseteq R$ denotes the (saturated)
homogeneous ideal of $X$ and $\mathcal I_X \subseteq \mathcal O_{\PP^3}$
denotes the ideal
sheaf of
$X$.
\item If $\Lambda$ is a $u \times (u+1)$ matrix with entries in a ring, we
denote by
$\Lambda_i$ ($i= 1, \dots, u+1$), the determinant of the submatrix of
$\Lambda$ obtained by deleting the $i-th$ column.
\end{itemize}
\medskip
Let $C \subseteq 2H$ be a curve of degree $d$ and genus $g$. Following
\cite{HS} we
associate to $C$ a triple $ Z \subset C'\subset P$ of subschemes of $H$,
where $C'$ and $P$ are (planar) curves and $Z$ is a
zero-dimensional subscheme; more precisely, $P$ is the pure $1$-dimensional
partof
$C\cap H$, $Z$
is the residual scheme to
$P$  in $C\cap H$, and $C'$ is the residual scheme to the intersection of
$C$ with $H$.
In other words we have an exact sequence:
\begin{equation}\label{residual sequence}
0\to \mathcal I_{C',\PP^3}(-1)\to \mathcal I_{C,\PP^3}\to \mathcal I_{C\cap
H,H}\to
0,
\end{equation}
where the first map is induced by multiplication by $x$ and $\mathcal
I_{C\cap H,H}= \mathcal I_{Z,H}(-P)$.
\medskip

We set:
\begin{itemize}
\item $I_{C'} = (x, p)$, where $p\in S$, and $\delta := \deg C'$
\item $I_P= (x, p h)$. Hence $d = 2\delta +\deg h$.
\item Let $A=(a_{i,j})$ be a fixed Hilbert-Burch matrix of $Z$ and let
$s\ge 1$ be
the number of rows of $A$.
Note that since $Z \subseteq C'$ there are forms
$p_1,\dots,p_{s+1} \in S$ (not uniquely determined) such that
$p = \sum_{i=1}^{s+1} p_iA_i$.
\item We have $I_{Z,H} = (A_1, \dots, A_{s+1})S$ and $I_{C\cap H,H} =p h
I_{Z,H}$.
\end{itemize}
Therefore the residual exact sequence (\ref{residual sequence}) can be
rewritten as \begin{equation}\label{ residual sequence2} 0\to \mathcal
I_{C',\PP^3}(-1)\to \mathcal I_{C,\PP^3}\to p hI_{Z,H}\to 0 \end{equation} and,
since $H^1_*(\mathcal I_C)= 0$, it is equivalent to the exact sequence
\begin{equation}\label{expected residual sequence}
0\to I_{C'}(-1)\xrightarrow{x} I_C \to p h I_{Z,H}\to 0 \end{equation}

\begin{remark} Notice that $Z\ne \emptyset$ if and only if $C$ is not aCM.
Furthermore if $Z=\emptyset$, we put $I_{Z,H}=S$ and the residual
sequence becomes
$$0\to (x,p)(-1)\xrightarrow{x} I_C \to p h S\to 0.$$

Then it is not difficult to see that the ideal of $C$ has the form
$I_C=(x^2,xp,ph+xf_1)$, i.e. it corresponds formally to a matrix
$B=[p,f_1]\in S^{1,2}$ as in Theorem \ref {thm-intro}, with $s=0$.
\end {remark}

>From now on we shall assume $Z\ne \emptyset$, i.e. that
$C$ is not aCM.

\bigskip

\section{Ideals with the expected residual sequence} \label{sec-resid-seq}
\bigskip
In this section we study ideals that behave, in some respects,  as the
ideal  of a curve
in a  double plane. We describe explicitly the minimal generators for them
as well as a
graded minimal free resolution. Now let us fix a triple $ Z \subset
C'\subset P$
of subschemes of a plane $H$, defined by the ideal $xR$, and use the above
notation.
\medskip
\begin{definition} Let $J \subseteq R$ be a homogeneous ideal.
We say that $J$ admits the {\it expected residual sequence} if it fits
into an exact sequence
\begin{equation}\label{expected residual sequence}
0\to I_{C'}(-1)\xrightarrow{x} J \to p h I_{Z,H}\to 0 \end{equation}
where the first map is multiplication by $x$.
\medskip
This means that $(J+xR)/xR = phI_{Z,H}$ and that $J:(x) = (x,p)$.
In other words $J$ behaves, with respect to the
data $p$, $h$ and $Z$, as the ideal of a curve in $2H$.
It is clear that $J \subseteq R$ has the expected residual sequence
if and only if there is
an exact sequence
\begin{equation}\label{expected residual sequence2}
0\to R/(x,p)(-1)\xrightarrow{x} R/J \to R /p hI_Z \to 0
\end{equation}
(compare with (\ref{residual sequence})).
\end{definition}

\begin{lemma} Let $J\subseteq R $ be a homogeneous ideal with the expected
residual
sequence. Then
$J = I_X$, where $X$ is a  one-dimensional closed subscheme of degree
$d$ and arithmetic genus
$\binom {d-\delta-1}{2}+\binom {\delta-1}{2}+\delta-\deg Z-1$
(i.e. $p_a(X)= p_a(P)+p_a(C')+\delta-\deg Z-1$).
\end{lemma}

\begin{proof} We show first that $J$ is saturated. Let $\ell$ be a general
element
of $R_1$ and let
$\phi \in R$ be such that $\ell\phi \in J$. Since $J+xR = phI_Z$ is
saturated,  it
follows that $\phi = \alpha + \beta x$, where $\alpha \in J$ and $\beta \in S$.
Then $\ell\beta x \in J$, whence $\ell\beta \in J:(x)$ which implies
$\beta \in J:(x)$, since $J:(x)$ is saturated.
Then $\beta x \in J$, whence $\phi \in J$ and $J$  is saturated.
To complete the proof it is sufficient to compute the Hilbert polynomial of
$X$
from the exact sequence (\ref{expected residual sequence2}). \end{proof}
\bigskip

Now we want to characterize the ideals with expected residual
sequence.

Let $A \in S^{s,s+1}$ be a Hilbert-Burch matrix
of $Z$.  Since $p \in I_{Z,H}$ we can fix
$p_1,\dots,p_{s+1} \in S$ such that
$$
p = \left |\begin{array}{c} p_1 \; \dots\; p_{s+1} \\A \end{array}\right|.
$$
For technical reasons we will write \begin{equation}\label{expression of p}
p = (-1)^s \left |\begin{array}{c} A\\ p_1\;\dots \;p_{s+1} \end{array}\right|
.\end{equation}
\medskip
\begin{remark} Let $J \subseteq R$ be an ideal having the expected residual
sequence.
Since $I_{Z,H}= (A_1, \dots, A_{s+1})S$ there are homogeneous polynomials
$G_1, \dots, G_{s+1} \in S$ such that
\begin{equation} \label{J}
J = (x^2, xp, p h A_1+xG_1, - p h A_2+xG_2,\; \dots\;, (-1)^s p h
A_{s+1}+xG_{s+1})
\end{equation}
and $\deg G_i = \deg (ph) + \deg A_i - 1$ for $i = 1, \dots, s+1$.

Notice that not every homogeneous ideal $J\subseteq R$ generated as above
has the
expected residual sequence. This depends on the choice of the $G_i$'s, as
we are going to
show.
\end{remark}
\medskip
\begin{proposition} \label{good residual} Consider a homogeneous ideal $J$
as in (\ref{J}).
Then the following conditions are equivalent:
\begin{itemize}
\item[(i)] $J$ has the expected residual sequence;
\item[(ii)] $J:(x) = (x,p)$;
\item[(iii)] there are homogeneous polynomials $f_1,\dots, f_s \in S$
such that
$$
A\left[\begin{array}{c}
G_1\\ \vdots \\ G_{s+1}
\end{array}
\right] = p \left [\begin{array}{c}
f_1\\ \vdots \\ f_s\\
\end{array}
\right]
$$
\item[(iv)]  there are homogeneous polynomials $f_1,\dots,f_{s+1}\in S$
such that the matrix
$$
B:= \left[\begin{array}{cc} A&\begin{array}{c}
f_1\\ \vdots
\end{array}\\  p_1\;\dots \;p_{s+1}&f_{s+1}
\end{array}
\right]
$$
is homogeneous and $G_i = (-1)^{i-1}B_i$ for $i = 1, \dots s+1$.
\end{itemize}
\end{proposition}

\begin{proof}
(i) $\Leftrightarrow$ (ii) is clear because $J(R/xR) = p h
I_{Z,H}$.

(ii) $\Rightarrow$ (iv). It is sufficient to show that there are
homogeneous polynomials $f_1,\dots,f_{s+1}\in S$ such that
\begin{equation}\label {b}
\left [\begin{array}{c} A\\ p_1\;\dots \;p_{s+1}
\end{array}\right]\left[\begin{array}{c}
G_1\\ \vdots \\ G_{s+1}\\
\end{array}
\right] = p \left [\begin{array}{c}
f_1\\ \vdots \\ f_{s+1}\\
\end{array}
\right]
\end{equation}
and then to apply Cramer's rule.

Now if $(\lambda_1, \dots \lambda_{s+1})$ is a row of
$A$ we have $\sum_{i=1}^{s+1}\lambda_i(-1)^{i-1}A_i=0$, whence
$$
x\sum_{i=1}^{s+1}\lambda_iG_i\in J
$$
which implies $\sum_{i=1}^{s+1}\lambda_iG_i \in (x,p)$. Since no
term of this sum contains $x$ we get $\sum_{i=1}^{s+1}\lambda_iG_i
\in pS$.

It remains to prove that $\sum_{i=1}^{s+1}p_iG_i \in pS$. As above
we have:

\begin{equation} \label{gamma2h}
p^2h + x\sum_{i=1}^{s+1}p_iG_i \in J
\end{equation}
whence it is sufficient to show that $p^2h \in J$.

Put $J' := J + p^2hR$. Then $J'+xR = J+xR$ by (\ref{gamma2h}).
Since $p^2h \in S$ we get $J':(x) = (x, p)$, and
hence there is an exact sequence (\ref{expected residual sequence})
with $J$ replaced by $J'$.
It follows that the inclusion map $J \hookrightarrow J'$ is bijective,
whence the conclusion.

(iv) $\Rightarrow$ (iii). By Cramer's rule we have that (i) implies (\ref
{b}), and
the conclusion is obvious.

(iii) $\Rightarrow$ (ii). It is sufficient to show that if $\phi \in S$ and
$x\phi \in J$
then $\phi \in (p)$.

An easy computation shows that there are
$\lambda_1, \dots \lambda_{s+1} \in S$ such that
$$
\phi = \sum_{i=1}^{s+1}\lambda_iG_i
$$
and
$$
\sum_{i=1}^{s+1}\lambda_i(-1)^{i-1}A_i = 0.
$$
whence the $(\lambda_1, \dots \lambda_{s+1})$ is a linear combination of
the rows of $A$.

By an easy computation using (iii) we get that $\phi$ is a multiple of $p$ and
(ii) follows.
\end{proof}

\begin{remark} \label{hp2} In the proof above we have seen:
If $J$ is an ideal with the expected residual sequence,
then $hp^2 \in J$.
\end{remark}
\medskip
\begin{corollary} \label{best generators} $J$ has the expected residual
sequence if and only if
$$
J = (x^2, xp, p h A_i+xB_i \mid i=1,\dots,s+1)
$$
where $B \in S^{s+1,s+2}$ is a matrix of the form
$$
B:= \left[\begin{array}{cc} A&\begin{array}{c}
f_1\\ \vdots
\end{array}\\ p_1\;\dots \;p_{s+1}&f_{s+1}
\end{array}
\right],
$$
with
$$
p = (-1)^s \left |\begin{array}{c} A\\ p_1\;\dots \;p_{s+1}
\end{array}\right|,$$
and
\begin {itemize}
\item $\deg f_i= \deg a_{i,j}+ \deg h + \deg A_j - 1$ for all $i=1, \dots ,
s$ and
for all $j = 1, \dots, s+1$,
\item $\deg f_{s+1}= \deg p_j + \deg h + \deg A_j - 1 = d - \delta -1$ for all
$j = 1, \dots, s+1$.
\end {itemize}
Moreover $\deg f_i \ge \deg a_{i,j} + \deg A_j - 1$ for all $i$ and $j$ as
above.

In particular $\deg f_i \ge \deg a_{i,j}$ for all $i,j$, where equality for
some
pair $(i,j)$ implies $\deg A_j = 1$ and $s=1$.
\end{corollary}
\begin{proof} Let $B$ be the matrix of Proposition \ref {good residual}. By
item
(iv) of the same Proposition we have $$ J = (x^2,
xp, p h A_1+xB_1, -p h A_2-xB_2,\; \dots\;, (-1)^{s+1}p h
A_{s+1}+(-1)^{s+1}xB_{s+1}).
$$
whence the first statement.
The computation for the degrees of $f_i$'s is an immediate consequence of
Proposition \ref {good residual}(iii).
The remaining statements follow by easy calculations.
\end{proof}
\medskip

Now we want to compute the minimal free resolution of a homogeneous ideal $J$
having the expected residual sequence. We fix some notation.

Set $J=(x^2,xp,\; phA_i+xB_i \mid 1\le i\leq s+1)$, where $B$ is a matrix
corresponding to $J$ (see Corollary \ref {best generators}) and let $M$ be
the matrix obtained from $B$ by deleting the last row, i.e.
$$
M:= \left[\begin{array}{cc} A&\begin{array}{c}
f_1\\ \vdots \\
f_{s}
\end{array}\\
\end{array}
\right] \in S^{s, s+2}.
$$

Recall that $A \in S^{s,s+1}$ is the Hilbert-Burch matrix of $Z$ and write
the minimal
free resolution of $I_{Z,H}$ over $S$ as
\begin{equation}\label{resolution of Z}
0\to \overline G\xrightarrow {^tA}\overline F\to I_{Z,H}\to 0. \end{equation}
Set $F=\overline F\otimes_S R$,
$G=\overline G\otimes_S R$.

Furthermore if $\Lambda$ is a given matrix and $t > 0$, we denote by
$I_t(\Lambda)$
the ideal generated by the minors of order $t$ of $\Lambda$.
Finally $E_t$ denotes the $t\times t$ identity matrix.
\bigskip
\begin{theorem}\label{resolution of J} Assume that $J$ has the expected
residual sequence, and let the notation be as above.
Then $J$ has the graded minimal free resolution
\begin{equation}
\begin{split}
0\to G(-d-1+\delta) \xrightarrow {\alpha_3}
&G(-d+\delta)\oplus F(-d-1+\delta)\oplus R(-2-\delta) \xrightarrow
{\alpha_2}\\ &R(-2)\oplus F(-d+\delta)\oplus R(-1-\delta)
\xrightarrow {\alpha_1}J\to 0 \end{split}
\end{equation}
\noindent where we have identified the maps with their matrices:
\bigskip
\flushleft $\begin{array}{lcrl}\alpha_1&=& \left [x^2,phA_1+xB_1,\right.
\dots,&(-1)^{i+1}\{phA_i+xB_i\},\dots ,\\
&&&\left.(-1)^s\{phA_{s+1}+xB_{s+1}\},
-px \right] \in R^{s+3}
\end{array}$

\bigskip

$\alpha_2=\left [\begin{array}{ccccccc}
0 & -B_1 & \dots&(-1)^iB_i & \dots&(-1)^{s+1}B_{s+1} & p\\ \\
^tM&&&xE_{s+1}&&&0\\ \\
& A_1h & \dots &(-1)^{i+1}A_ih & \dots & (-1)^s A_{s+1}h&x
\end{array}\right ]\in R^{s+3,2s+2}$

\bigskip

$\alpha_3=\left [\begin{array}{c} -xE_s\\ \\
^tM
\end{array}\right ]\in R^{2s+2,s}.$
\medskip
\end{theorem}
\begin{proof} Developing determinants along a row and using
$p=(-1)^{s+2}B_{s+2}$
it is easy to check that the sequence described above is a complex.
To check that it is exact we use the Buchsbaum-Eisenbud criterion (see
\cite {Ei},
Theorem 20.9). Since $I_s(\alpha_3)\supseteq x^sR+I_s(A)=x^sR+I_{Z,H}R$,
this ideal
has at least codimension $3$.
It remains to show that $I_{s+2}(\alpha_2)$ has at least codimension $2$.
It is immediate to see that $x^{s+2}\in I_{s+2}(\alpha_2)$,
that is $x \in \rad(I_{s+2}(\alpha_2))$.
Moreover if $N$ is the matrix obtained from $^tM$ by
deleting the first row we have
$$
\det \left [\begin{array}{cccc}
0&&(-1)^{s+1}B_{s+1}&p\\
&&0&0\\
N&&\vdots&\vdots\\
&&0&0\\
&&x&0\\
&&(-1)^sA_{s+1}h&x
\end{array}\right ]\in I_{s+2}(\alpha_2),
$$
\bigskip
\noindent whence $0 \neq A_1A_{s+1}hp \in \rad(I_{s+2}(\alpha_2))$.
This proves that $\codim (I_{s+2}(\alpha_2))\ge 2$.\par
Since the non-zero entries in the matrices $\alpha_1,\alpha_2,\alpha_3$
have positive degree,
the exact sequence above is a minimal free resolution.
\end{proof}

\bigskip

\section{Curves in the double plane} \label{curves}

\bigskip

In this section we relate ideals with expected residual sequence to curves in
a double plane. We use the notation as in Theorem \ref {resolution of J}.

Moreover for a graded $R$-module $M$ we denote by $M^* := \Hom_R(M,R)$
the $R$-dual of $M$ and by $M^\vee : = \Hom_K(M,K)$ the graded $K$-dual of
$M$.

The main result in this section is:

\begin{theorem}\label{curve and minors of M} Let $C$ be the scheme defined
by an
ideal $J$ having expected residual sequence. Then we have:
\begin{itemize}
\item[(i)] $H^1_*{(\mathcal I_C)}^\vee \cong \coker (\alpha_3^*)(-4)$;
\item[(ii)] $C$ is a curve if and only if the maximal minors of $M$ do not have
common zeroes in $H$.
\item[(iii)] If $C$ is a curve, then its Hartshorne-Rao
module $M_C := H^1_*{(\mathcal I_C)}$ is:
$$
M_C\cong \coker (G^*(-\delta-2)\oplus F^*(-\delta-1)\oplus
R(-d+\delta)\xrightarrow {\left [\begin{array}{cc} xE_s&M
\end{array}\right ]} G^*(-\delta-1)).
$$
\end{itemize}
\end{theorem}

\begin{proof} By duality $H^1_*{(\mathcal I_C)}^\vee \cong
\Ext^2_R(J,R)(-4)$ and
(i) follows by  Theorem \ref{resolution of J}. Moreover $C$ is a curve if
and only if
the graded $R$-module $H^1_*(\mathcal I_C)$ has finite length, whence
(ii) follows
from (i). Finally, since $M_C^{\vee} \cong M_C(d-2)$ (see \cite{HS},
Cor. 6.2) (iii) follows from (i).
\end{proof}

\begin{corollary}\label{mingen} $M_C$ is minimally
generated by $s$ homogeneous elements. Moreover if $s=1$, then $M_C$ is a
Koszul module, i.e.\ it has only $4$ relations.
\end{corollary}

\begin{remark} The self-duality of $M_C$, proved in \cite {HS} and used in
the proof of
Theorem \ref{curve and minors of M} could also be proved easily as a
consequence of
Theorem \ref{resolution of J}.
\end{remark}

The previous Theorem and Corollary \ref{best generators} imply part (a) of
Theorem \ref{thm-intro} of the introduction.
Now we are going to prove part (b) of Theorem \ref{thm-intro} and Corollary
\ref{cor-intro}.
As a byproduct we find an alternative proof of a key result in \cite {HS}.

\medskip

\begin{lemma}\label{lemma on M} Let $A \in S^{s,s+1}$ be a homogeneous
matrix. Then
the following conditions are equivalent:
\begin{itemize}
\item[(i)] there is a homogeneous matrix $M$ of the form
$$
M=\left [\begin{array}{ccccc}&&&f_1\\
A&&&&\\
&&&f_s
\end{array}\right] \in S^{s,s+2}
$$
such that $I_s(M)$ is irrelevant;
\item[(ii)] the subscheme $Z \subseteq \Proj(S)\cong \PP^2$ defined
by
$I_s(A)$ is zero-dimensional and locally a complete intersection;
\end{itemize}
\end{lemma}

\begin{proof} (i) $\Rightarrow$ (ii). Since $I_s(M)$ defines a subscheme of
codimension 3,  by \cite{B} $Z$ has codimension 2, hence is zero-dimensional.
Moreover, since $I_{s-1}(A)\supseteq I_s(M),  Z$ is locally a
complete intersection by \cite{KMMNP}, Proposition 3.2.

(ii) $\Rightarrow$ (i). By \cite {KMNP} (see also \cite{KMMNP}, Proposition
3.2.)
we may assume that the
$(s-1)\times (s-1)$ minors $A'_1, \dots, A'_{s+1}$ corresponding to the first
$s-1$ rows of $A$ generate an irrelevant ideal. Now consider the matrix
$$
M:= \left [\begin{array}{ccccc} &&&0\\
A&&&\vdots\\
&&&0\\
&&&f_s
\end{array}\right ],
$$
where $f_s \in S$ is a homogeneous polynomial  not vanishing at any point
of $Z$.
Then $I_s(M)\supseteq I_{Z,H}+ f_s(A'_1, \dots, A'_{s+1})$, and this
easily implies that $I_s(M)$ is irrelevant.
\end{proof}

\begin{corollary}\label{locally CI} {\rm (\cite {HS})} Let
$Z\subset C'\subset P$ be a triple as in Section 2 and let
$A$ be a Hilbert-Burch matrix of $Z$. Then the following conditions are
equivalent:
\begin{itemize}
\item[(i)] there is a curve $C \subseteq 2H$ corresponding to the triple;
\item[(ii)] $Z$ is locally a complete intersection.
\end{itemize}
\end{corollary}

\begin{proof} It follows easily from Corollary \ref{best generators}, Theorem
\ref{curve and minors of M}, and Lemma \ref{lemma on M}.
\end{proof}

\begin{corollary}\label{presentation} Let $A \in S^{s,s+1}$ and
$M := \left [\begin{array}{ccccc} &&&f_1\\
A&&&\vdots\\
&&&f_s
\end{array}\right ] \in S^{s+1,s+2}$
be homogeneous matrices such that
$\deg f_1 \ge \deg a_{1,1} + \sum_{j=1}^s \deg a_{j,j+1}-1$ and $I_s(M)$ is
irrelevant.

Then there is a curve in $2H$ whose Hartshorne-Rao module is presented by
the matrix
$\left [xE_s\; M \right]$.
\end{corollary}

\begin{proof} By Lemma \ref{lemma on M} $A$ is the Hilbert-Burch
matrix of a zero-dimensional scheme $Z \subseteq H$. Fix a
homogeneous polynomial $h \in S$ such that
$\deg h = \deg f_1 - (\deg a_{1,1} + \sum_{j=1}^s \deg a_{j,j+1}-1)$ and
construct a homogeneous matrix
$$
B:= \left[\begin{array}{cc} A&\begin{array}{c}
f_1\\ \vdots
\end{array}\\ p_1\;\dots \; p_{s+1}&f_{s+1}
\end{array}
\right],
$$
such that $p:= (-1)^s B_{s+2}$ is non-zero. Then the homogeneous ideal
$$
J = (x^2, xp, p h A_i+xB_i \mid i=1,\dots,s+1)
$$
defines a curve with the required properties by Lemma \ref{lemma on M} and
Theorem
\ref{curve and minors of M}.
\end{proof}

The above Corollary completes the proof of Theorem \ref{thm-intro}
(b).

Hence Theorem \ref{thm-intro} is proved. Moreover Corollary
\ref{cor-intro} follows by Theorem \ref{curve and minors of M} and
Corollary \ref {presentation}.

\begin{remark}\label {row} \begin {itemize}
\item[(i)] Let $C$ and $B$ be as in Corollary \ref{best
generators}. It is easy to see that elementary row operations on
$B$ produce the same curve, except adding to one of the first $s$ rows a
multiple of
the last row: such an operation can change $Z$, hence the curve.
\item[(ii)] Similarly elementary column operations on $B$
produce the same curve, unless $s=1$ and $\deg f_1 = \deg a_{1,1}$ or $\deg
f_1 = \deg
a_{1,2}$. Indeed by the degree conditions of Corollary \ref{best
generators} it is not
possible to add to one of the first $s+1$ columns a multiple of the last
one, except in
the particular case just mentioned. Any other column operation does not
change the
ideal.
\end{itemize}
\end{remark}

\bigskip

We end this section by proving Proposition \ref{prop-intro}.

We fix a non-aCM curve $C \subseteq 2H$. This
also fixes a presentation of its Hartshorne-Rao module $M_C$ as shown in
Theorem \ref{curve and minors of M}: observe that $M_C$
has a natural structure as $S$-module, and as such is represented by the
matrix $M$.
First we observe:

\begin{lemma} \label{ann M_C} We have:
\begin{itemize}
\item [(i)] $\Ann_S(M_C) = I_s(M)$;
\item [(ii)] $\Ann_R(M_C) = xR + I_s(M)R$.
\end{itemize}
\end{lemma}

\begin{proof} Both claims follow immediately by \cite{BE}.
\end{proof}

\begin{proof}[Proof of Proposition \ref {prop-intro}]
\mbox{}

(i) $\Leftrightarrow$ (iii) is a consequence of Corollary \ref{cor-intro} and
Corollary
\ref{presentation}.

(ii) $\Rightarrow$ (i) follows from well-known facts on extremal curves
(see e.g. \cite
{MDP}).

(i) $\Rightarrow$ (ii). Let $C \in \cL_N$, $C \subseteq 2 H$ and let $D \in
\cL_N$ be
minimal.
 Then $D$ must
lie on some quadric $Q$ (see, e.g., \cite {MDPA}, Ch. III, Prop. 3.6 and
Th. 5.1). We
are going to show that if $C$ is not extremal, then $Q$ is exactly $2H$.

Set $qR =
I_Q$. Then it is easy to see that $q \in \Ann_R(M_D) =\Ann_R(M_C)$.

Assume first that $Q$ is irreducible. Then $Q$ must be smooth, because $D$ is
not aCM. Then $D$ is of type $(0,b)$ (say) where $b := \deg D$. It  follows
that $M_D$ is minimally generated by $b-1$ elements in degree zero (see,
e.g., \cite
{GM}). This implies, by Corollary \ref{mingen}, that
$b-1=s$. Now $q \in \Ann_R(M_C)$ is an irreducible form of degree $2$,
whence $s
\le 2$ by Lemma \ref {ann M_C}. If $s=1$ then $D$ is a curve of degree $2$ and
genus $-1$, hence extremal. Thus, we get $s = 2$ and
the Rao function of $D$ is easily seen to be: $\rho(0) = \rho(1) =2$
and $\rho(j) = 0$ for $j \not = 0,1$ (see, e.g., \cite {GM}).
On the other hand by Theorem \ref{curve and minors of M} and the degree
conditions of
Corollary \ref {best generators} it follows that there is at least one $j$
such that
$\rho_C(j)\ge 3$, a contradiction.

Assume now that $Q$ is the union of two
distinct planes. Then $D$ is extremal by \cite {H}, Example 5.7 or \cite {HS},
Prop.\ 9.5.

Therefore $Q$ is a double plane $2H'$. Since  $I_H$ and $I_{H'}$ are
contained in
$\Ann_R (M_D)$  by Lemma \ref {ann M_C} and  since $D$ is not extremal, we have
$H = H'$.
\end{proof}

\end{document}